\catcode`\@=11
\font\tensmc=cmcsc10      
\def\smc{\tensmc}

\def\hcorrection#1{\advance\hoffset by #1 }
\def\vcorrection#1{\advance\voffset by #1 }
\def\wlog#1{}
\newif\iftitle@
\outer\def\title{\title@true\vglue 24\p@ plus 12\p@ minus 12\p@
   \bgroup\let\\=\cr\tabskip\centering
   \halign to \hsize\bgroup\tenbf\hfill\ignorespaces##\unskip\hfill\cr}
\def\endtitle{\cr\egroup\egroup\vglue 18\p@ plus 12\p@ minus 6\p@}
\outer\def\author{\iftitle@\vglue -18\p@ plus -12\p@ minus -6\p@\fi\vglue
    12\p@ plus 6\p@ minus 3\p@\bgroup\let\\=\cr\tabskip\centering
    \halign to \hsize\bgroup\smc\hfill\ignorespaces##\unskip\hfill\cr}
\def\endauthor{\cr\egroup\egroup\vglue 18\p@ plus 12\p@ minus 6\p@}
\outer\def\heading{\bigbreak\bgroup\let\\=\cr\tabskip\centering
    \halign to \hsize\bgroup\smc\hfill\ignorespaces##\unskip\hfill\cr}
\def\endheading{\cr\egroup\egroup\nobreak\medskip}

\outer\def\proclaim#1{\medbreak\noindent\smc\ignorespaces
    #1\unskip.\enspace\sl\ignorespaces}
\outer\def\endproclaim{\par\ifdim\lastskip<\medskipamount\removelastskip
  \penalty 55 \fi\medskip\rm}
\outer\def\demo#1{\par\ifdim\lastskip<\smallskipamount\removelastskip
    \smallskip\fi\noindent{\smc\ignorespaces#1\unskip:\enspace}\rm
      \ignorespaces}

\newcount\footmarkcount@
\footmarkcount@=1
\def\makefootnote@#1#2{\insert\footins{\interlinepenalty=100
  \splittopskip=\ht\strutbox \splitmaxdepth=\dp\strutbox
  \floatingpenalty=\@MM
  \leftskip=\z@\rightskip=\z@\spaceskip=\z@\xspaceskip=\z@
  \noindent{#1}\footstrut\rm\ignorespaces #2\strut}}
\def\footnote{\let\@sf=\empty\ifhmode\edef\@sf{\spacefactor
   =\the\spacefactor}\/\fi\futurelet\next\footnote@}
\def\footnote@{\ifx"\next\let\next\footnote@@\else
    \let\next\footnote@@@\fi\next}
\def\footnote@@"#1"#2{#1\@sf\relax\makefootnote@{#1}{#2}}
\def\footnote@@@#1{$^{\number\footmarkcount@}$\makefootnote@
   {$^{\number\footmarkcount@}$}{#1}\global\advance\footmarkcount@ by 1 }

\hyphenation{man-u-script man-u-scripts ap-pen-dix ap-pen-di-ces}
\hyphenation{data-base data-bases}
\ifx\amstexloaded@\relax\catcode`\@=13
  \endinput\else\let\amstexloaded@=\relax\fi
\newlinechar=`\^^J
\def\eat@#1{}
\def\Space@.{\futurelet\Space@\relax}
\Space@. %
\newhelp\athelp@
{Only certain combinations beginning with @ make sense to me.^^J
Perhaps you wanted \string\@\space for a printed @?^^J
I've ignored the character or group after @.}
\def\futureletnextat@{\futurelet\next\at@}
{\catcode`\@=\active
\lccode`\Z=`\@ \lowercase
{\gdef@{\expandafter\csname futureletnextatZ\endcsname}
\expandafter\gdef\csname atZ\endcsname
   {\ifcat\noexpand\next a\def\next{\csname atZZ\endcsname}\else
   \ifcat\noexpand\next0\def\next{\csname atZZ\endcsname}\else
    \def\next{\csname atZZZ\endcsname}\fi\fi\next}
\expandafter\gdef\csname atZZ\endcsname#1{\expandafter
   \ifx\csname #1Zat\endcsname\relax\def\next
     {\errhelp\expandafter=\csname athelpZ\endcsname
      \errmessage{Invalid use of \string@}}\else
       \def\next{\csname #1Zat\endcsname}\fi\next}
\expandafter\gdef\csname atZZZ\endcsname#1{\errhelp
    \expandafter=\csname athelpZ\endcsname
      \errmessage{Invalid use of \string@}}}}
\def\atdef@#1{\expandafter\def\csname #1@at\endcsname}
\newhelp\defahelp@{If you typed \string\define\space cs instead of
\string\define\string\cs\space^^J
I've substituted an inaccessible control sequence so that your^^J
definition will be completed without mixing me up too badly.^^J
If you typed \string\define{\string\cs} the inaccessible control sequence^^J
was defined to be \string\cs, and the rest of your^^J
definition appears as input.}
\newhelp\defbhelp@{I've ignored your definition, because it might^^J
conflict with other uses that are important to me.}
\def\define{\futurelet\next\define@}
\def\define@{\ifcat\noexpand\next\relax
  \def\next{\define@@}%
  \else\errhelp=\defahelp@
  \errmessage{\string\define\space must be followed by a control
     sequence}\def\next{\def\garbage@}\fi\next}
\def\undefined@{}
\def\preloaded@{}
\def\define@@#1{\ifx#1\relax\errhelp=\defbhelp@
   \errmessage{\string#1\space is already defined}\def\next{\def\garbage@}%
   \else\expandafter\ifx\csname\expandafter\eat@\string
         #1@\endcsname\undefined@\errhelp=\defbhelp@
   \errmessage{\string#1\space can't be defined}\def\next{\def\garbage@}%
   \else\expandafter\ifx\csname\expandafter\eat@\string#1\endcsname\relax
     \def\next{\def#1}\else\errhelp=\defbhelp@
     \errmessage{\string#1\space is already defined}\def\next{\def\garbage@}%
      \fi\fi\fi\next}
\def\famzero{\fam\z@}

\def\lim{\mathop{\famzero lim}}

\def\textfont@#1#2{\def#1{\relax\ifmmode
    \errmessage{Use \string#1\space only in text}\else#2\fi}}
\textfont@\rm\tenrm
\textfont@\it\tenit
\textfont@\sl\tensl
\textfont@\bf\tenbf
\textfont@\smc\tensmc
\let\ic@=\/
\def\/{\unskip\ic@}
\def\textfonti{\the\textfont1 }
\def\t#1#2{{\edef\next{\the\font}\textfonti\accent"7F \next#1#2}}
\let\B=\=
\let\D=\.
\def~{\unskip\nobreak\ \ignorespaces}
{\catcode`\@=\active
\gdef\@{\char'100 }}
\atdef@-{\leavevmode\futurelet\next\athyph@}
\def\athyph@{\ifx\next-\let\next=\athyph@@
  \else\let\next=\athyph@@@\fi\next}
\def\athyph@@@{\hbox{-}}
\def\athyph@@#1{\futurelet\next\athyph@@@@}
\def\athyph@@@@{\if\next-\def\next##1{\hbox{---}}\else
    \def\next{\hbox{--}}\fi\next}
\def\.{.\spacefactor=\@m}
\atdef@.{\null.}
\atdef@,{\null,}
\atdef@;{\null;}
\atdef@:{\null:}
\atdef@?{\null?}
\atdef@!{\null!}
\def\srdr@{\thinspace}
\def\drsr@{\kern.02778em}
\def\sldl@{\kern.02778em}
\def\dlsl@{\thinspace}
\atdef@"{\unskip\futurelet\next\atqq@}
\def\atqq@{\ifx\next\Space@\def\next. {\atqq@@}\else
         \def\next.{\atqq@@}\fi\next.}
\def\atqq@@{\futurelet\next\atqq@@@}
\def\atqq@@@{\ifx\next`\def\next`{\atqql@}\else\def\next'{\atqqr@}\fi\next}
\def\atqql@{\futurelet\next\atqql@@}
\def\atqql@@{\ifx\next`\def\next`{\sldl@``}\else\def\next{\dlsl@`}\fi\next}
\def\atqqr@{\futurelet\next\atqqr@@}
\def\atqqr@@{\ifx\next'\def\next'{\srdr@''}\else\def\next{\drsr@'}\fi\next}

\def\textfontii{\the\textfont2 }
\def\{{\relax\ifmmode\lbrace\else
    {\textfontii f}\spacefactor=\@m\fi}
\def\}{\relax\ifmmode\rbrace\else
    \let\@sf=\empty\ifhmode\edef\@sf{\spacefactor=\the\spacefactor}\fi
      {\textfontii g}\@sf\relax\fi}
\def\nonhmodeerr@#1{\errmessage
     {\string#1\space allowed only within text}}
\def\linebreak{\relax\ifhmode\unskip\break\else
    \nonhmodeerr@\linebreak\fi}
\def\allowlinebreak{\relax
   \ifhmode\allowbreak\else\nonhmodeerr@\allowlinebreak\fi}
\newskip\saveskip@
\def\nolinebreak{\relax\ifhmode\saveskip@=\lastskip\unskip
  \nobreak\ifdim\saveskip@>\z@\hskip\saveskip@\fi
   \else\nonhmodeerr@\nolinebreak\fi}
\def\newline{\relax\ifhmode\null\hfil\break
    \else\nonhmodeerr@\newline\fi}
\def\nonmathaerr@#1{\errmessage
     {\string#1\space is not allowed in display math mode}}
\def\nonmathberr@#1{\errmessage{\string#1\space is allowed only in math mode}}
\def\mathbreak{\relax\ifmmode\ifinner\break\else
   \nonmathaerr@\mathbreak\fi\else\nonmathberr@\mathbreak\fi}
\def\nomathbreak{\relax\ifmmode\ifinner\nobreak\else
    \nonmathaerr@\nomathbreak\fi\else\nonmathberr@\nomathbreak\fi}
\def\allowmathbreak{\relax\ifmmode\ifinner\allowbreak\else
     \nonmathaerr@\allowmathbreak\fi\else\nonmathberr@\allowmathbreak\fi}
\def\pagebreak{\relax\ifmmode
   \ifinner\errmessage{\string\pagebreak\space
     not allowed in non-display math mode}\else\postdisplaypenalty-\@M\fi
   \else\ifvmode\penalty-\@M\else\edef\spacefactor@
       {\spacefactor=\the\spacefactor}\vadjust{\penalty-\@M}\spacefactor@
        \relax\fi\fi}
\def\nopagebreak{\relax\ifmmode
     \ifinner\errmessage{\string\nopagebreak\space
    not allowed in non-display math mode}\else\postdisplaypenalty\@M\fi
    \else\ifvmode\nobreak\else\edef\spacefactor@
        {\spacefactor=\the\spacefactor}\vadjust{\penalty\@M}\spacefactor@
         \relax\fi\fi}
\def\newpage{\relax\ifvmode\vfill\penalty-\@M\else\nonvmodeerr@\newpage\fi}
\def\nonvmodeerr@#1{\errmessage
    {\string#1\space is allowed only between paragraphs}}
\def\smallpagebreak{\relax\ifvmode\smallbreak
      \else\nonvmodeerr@\smallpagebreak\fi}
\def\medpagebreak{\relax\ifvmode\medbreak
       \else\nonvmodeerr@\medpagebreak\fi}
\def\bigpagebreak{\relax\ifvmode\bigbreak
      \else\nonvmodeerr@\bigpagebreak\fi}
\newdimen\captionwidth@
\captionwidth@=\hsize
\advance\captionwidth@ by -1.5in
\def\caption#1{}
\def\topspace#1{\gdef\thespace@{#1}\ifvmode\def\next
    {\futurelet\next\topspace@}\else\def\next{\nonvmodeerr@\topspace}\fi\next}
\def\topspace@{\ifx\next\Space@\def\next. {\futurelet\next\topspace@@}\else
     \def\next.{\futurelet\next\topspace@@}\fi\next.}
\def\topspace@@{\ifx\next\caption\let\next\topspace@@@\else
    \let\next\topspace@@@@\fi\next}
 \def\topspace@@@@{\topinsert\vbox to
       \thespace@{}\endinsert}
\def\topspace@@@\caption#1{\topinsert\vbox to
    \thespace@{}\nobreak
      \smallskip
    \setbox\z@=\hbox{\noindent\ignorespaces#1\unskip}%
   \ifdim\wd\z@>\captionwidth@
   \centerline{\vbox{\hsize=\captionwidth@\noindent\ignorespaces#1\unskip}}%
   \else\centerline{\box\z@}\fi\endinsert}
\def\midspace#1{\gdef\thespace@{#1}\ifvmode\def\next
    {\futurelet\next\midspace@}\else\def\next{\nonvmodeerr@\midspace}\fi\next}
\def\midspace@{\ifx\next\Space@\def\next. {\futurelet\next\midspace@@}\else
     \def\next.{\futurelet\next\midspace@@}\fi\next.}
\def\midspace@@{\ifx\next\caption\let\next\midspace@@@\else
    \let\next\midspace@@@@\fi\next}
 \def\midspace@@@@{\midinsert\vbox to
       \thespace@{}\endinsert}
\def\midspace@@@\caption#1{\midinsert\vbox to
    \thespace@{}\nobreak
      \smallskip
      \setbox\z@=\hbox{\noindent\ignorespaces#1\unskip}%
      \ifdim\wd\z@>\captionwidth@
    \centerline{\vbox{\hsize=\captionwidth@\noindent\ignorespaces#1\unskip}}%
    \else\centerline{\box\z@}\fi\endinsert}
\mathchardef\prime@="0230
\def\prime{{{}\prime@{}}}
\def\prim@s{\prime@\futurelet\next\pr@m@s}

\def\,{\relax\ifmmode\mskip\thinmuskip\else\thinspace\fi}
\def\!{\relax\ifmmode\mskip-\thinmuskip\else\negthinspace\fi}
\def\frac#1#2{{#1\over#2}}

\def\:{\nobreak\hskip.1111em{:}\hskip.3333em plus .0555em\relax}
\def\intic@{\mathchoice{\hskip5\p@}{\hskip4\p@}{\hskip4\p@}{\hskip4\p@}}
\def\negintic@
 {\mathchoice{\hskip-5\p@}{\hskip-4\p@}{\hskip-4\p@}{\hskip-4\p@}}
\def\intkern@{\mathchoice{\!\!\!}{\!\!}{\!\!}{\!\!}}
\def\intdots@{\mathchoice{\cdots}{{\cdotp}\mkern1.5mu
    {\cdotp}\mkern1.5mu{\cdotp}}{{\cdotp}\mkern1mu{\cdotp}\mkern1mu
      {\cdotp}}{{\cdotp}\mkern1mu{\cdotp}\mkern1mu{\cdotp}}}
\newcount\intno@
\def\iint{\intno@=\tw@\futurelet\next\ints@}
\def\iiint{\intno@=\thr@@\futurelet\next\ints@}
\def\iiiint{\intno@=4 \futurelet\next\ints@}
\def\idotsint{\intno@=\z@\futurelet\next\ints@}
\def\ints@{\findlimits@\ints@@}
\newif\iflimtoken@
\newif\iflimits@
\def\findlimits@{\limtoken@false\limits@false\ifx\next\limits
 \limtoken@true\limits@true\else\ifx\next\nolimits\limtoken@true\limits@false
    \fi\fi}
\def\multintlimits@{\intop\ifnum\intno@=\z@\intdots@
  \else\intkern@\fi
    \ifnum\intno@>\tw@\intop\intkern@\fi
     \ifnum\intno@>\thr@@\intop\intkern@\fi\intop}
\def\multint@{\int\ifnum\intno@=\z@\intdots@\else\intkern@\fi
   \ifnum\intno@>\tw@\int\intkern@\fi
    \ifnum\intno@>\thr@@\int\intkern@\fi\int}
\def\ints@@{\iflimtoken@\def\ints@@@{\iflimits@
   \negintic@\mathop{\intic@\multintlimits@}\limits\else
    \multint@\nolimits\fi\eat@}\else
     \def\ints@@@{\multint@\nolimits}\fi\ints@@@}
\def\Sb{_\bgroup\vspace@
        \baselineskip=\fontdimen10 \scriptfont\tw@
        \advance\baselineskip by \fontdimen12 \scriptfont\tw@
        \lineskip=\thr@@\fontdimen8 \scriptfont\thr@@
        \lineskiplimit=\thr@@\fontdimen8 \scriptfont\thr@@
        \Let@\vbox\bgroup\halign\bgroup \hfil$\scriptstyle
            {##}$\hfil\cr}
\def\endSb{\crcr\egroup\egroup\egroup}
\def\Sp{^\bgroup\vspace@
        \baselineskip=\fontdimen10 \scriptfont\tw@
        \advance\baselineskip by \fontdimen12 \scriptfont\tw@
        \lineskip=\thr@@\fontdimen8 \scriptfont\thr@@
        \lineskiplimit=\thr@@\fontdimen8 \scriptfont\thr@@
        \Let@\vbox\bgroup\halign\bgroup \hfil$\scriptstyle
            {##}$\hfil\cr}
\def\endSp{\crcr\egroup\egroup\egroup}
\def\Let@{\relax\iffalse{\fi\let\\=\cr\iffalse}\fi}
\def\vspace@{\def\vspace##1{\noalign{\vskip##1 }}}
\def\aligned{\,\vcenter\bgroup\vspace@\Let@\openup\jot\m@th\ialign
  \bgroup \strut\hfil$\displaystyle{##}$&$\displaystyle{{}##}$\hfil\crcr}
\def\endaligned{\crcr\egroup\egroup}
\def\matrix{\,\vcenter\bgroup\Let@\vspace@
    \normalbaselines
  \m@th\ialign\bgroup\hfil$##$\hfil&&\quad\hfil$##$\hfil\crcr
    \mathstrut\crcr\noalign{\kern-\baselineskip}}
\def\endmatrix{\crcr\mathstrut\crcr\noalign{\kern-\baselineskip}\egroup
                \egroup\,}
\newtoks\hashtoks@
\hashtoks@={#}
\def\format{\crcr\egroup\iffalse{\fi\ifnum`}=0 \fi\format@}
\def\format@#1\\{\def\preamble@{#1}%
  \def\c{\hfil$\the\hashtoks@$\hfil}%
  \def\r{\hfil$\the\hashtoks@$}%
  \def\l{$\the\hashtoks@$\hfil}%
  \setbox\z@=\hbox{\xdef\Preamble@{\preamble@}}\ifnum`{=0 \fi\iffalse}\fi
   \ialign\bgroup\span\Preamble@\crcr}
 
\let\hdots=\ldots
\def\cases{\left\{\,\vcenter\bgroup\vspace@
     \normalbaselines\openup\jot\m@th
       \Let@\ialign\bgroup$##$\hfil&\quad$##$\hfil\crcr
      \mathstrut\crcr\noalign{\kern-\baselineskip}}

\newif\iftagsleft@
\tagsleft@true
\def\TagsOnRight{\global\tagsleft@false}
\def\tag#1$${\iftagsleft@\leqno\else\eqno\fi
 \hbox{\def\pagebreak{\global\postdisplaypenalty-\@M}%
 \def\nopagebreak{\global\postdisplaypenalty\@M}\rm(#1\unskip)}%
  $$\postdisplaypenalty\z@\ignorespaces}
\interdisplaylinepenalty=\@M
\def\allowdisplaybreak@{\def\allowdisplaybreak{\noalign{\allowbreak}}}
\def\displaybreak@{\def\displaybreak{\noalign{\break}}}
\def\align#1\endalign{\def\tag{&}\vspace@\allowdisplaybreak@\displaybreak@
  \iftagsleft@\lalign@#1\endalign\else
   \ralign@#1\endalign\fi}
\def\ralign@#1\endalign{\displ@y\Let@\tabskip\centering\halign to\displaywidth
     {\hfil$\displaystyle{##}$\tabskip=\z@&$\displaystyle{{}##}$\hfil
       \tabskip=\centering&\llap{\hbox{(\rm##\unskip)}}\tabskip\z@\crcr
             #1\crcr}}
\def\lalign@
 #1\endalign{\displ@y\Let@\tabskip\centering\halign to \displaywidth
   {\hfil$\displaystyle{##}$\tabskip=\z@&$\displaystyle{{}##}$\hfil
   \tabskip=\centering&\kern-\displaywidth
        \rlap{\hbox{(\rm##\unskip)}}\tabskip=\displaywidth\crcr
               #1\crcr}}
\def\overrightarrow{\mathpalette\overrightarrow@}
\def\overrightarrow@#1#2{\vbox{\ialign{$##$\cr
    #1{-}\mkern-6mu\cleaders\hbox{$#1\mkern-2mu{-}\mkern-2mu$}\hfill
     \mkern-6mu{\to}\cr
     \noalign{\kern -1\p@\nointerlineskip}
     \hfil#1#2\hfil\cr}}}
\def\overleftarrow{\mathpalette\overleftarrow@}
\def\overleftarrow@#1#2{\vbox{\ialign{$##$\cr
     #1{\leftarrow}\mkern-6mu\cleaders\hbox{$#1\mkern-2mu{-}\mkern-2mu$}\hfill
      \mkern-6mu{-}\cr
     \noalign{\kern -1\p@\nointerlineskip}
     \hfil#1#2\hfil\cr}}}
\def\overleftrightarrow{\mathpalette\overleftrightarrow@}
\def\overleftrightarrow@#1#2{\vbox{\ialign{$##$\cr
     #1{\leftarrow}\mkern-6mu\cleaders\hbox{$#1\mkern-2mu{-}\mkern-2mu$}\hfill
       \mkern-6mu{\to}\cr
    \noalign{\kern -1\p@\nointerlineskip}
      \hfil#1#2\hfil\cr}}}
\def\underrightarrow{\mathpalette\underrightarrow@}
\def\underrightarrow@#1#2{\vtop{\ialign{$##$\cr
    \hfil#1#2\hfil\cr
     \noalign{\kern -1\p@\nointerlineskip}
    #1{-}\mkern-6mu\cleaders\hbox{$#1\mkern-2mu{-}\mkern-2mu$}\hfill
     \mkern-6mu{\to}\cr}}}
\def\underleftarrow{\mathpalette\underleftarrow@}
\def\underleftarrow@#1#2{\vtop{\ialign{$##$\cr
     \hfil#1#2\hfil\cr
     \noalign{\kern -1\p@\nointerlineskip}
     #1{\leftarrow}\mkern-6mu\cleaders\hbox{$#1\mkern-2mu{-}\mkern-2mu$}\hfill
      \mkern-6mu{-}\cr}}}
\def\underleftrightarrow{\mathpalette\underleftrightarrow@}
\def\underleftrightarrow@#1#2{\vtop{\ialign{$##$\cr
      \hfil#1#2\hfil\cr
    \noalign{\kern -1\p@\nointerlineskip}
     #1{\leftarrow}\mkern-6mu\cleaders\hbox{$#1\mkern-2mu{-}\mkern-2mu$}\hfill
       \mkern-6mu{\to}\cr}}}
\def\sqrt#1{\radical"270370 {#1}}
\def\dots{\relax\ifmmode\let\next=\ldots\else\let\next=\tdots@\fi\next}
\def\tdots@{\unskip\ \tdots@@}
\def\tdots@@{\futurelet\next\tdots@@@}
\def\tdots@@@{$\mathinner{\ldotp\ldotp\ldotp}\,
   \ifx\next,$\else
   \ifx\next.\,$\else
   \ifx\next;\,$\else
   \ifx\next:\,$\else
   \ifx\next?\,$\else
   \ifx\next!\,$\else
   $ \fi\fi\fi\fi\fi\fi}
\def\text{\relax\ifmmode\let\next=\text@\else\let\next=\text@@\fi\next}
\def\text@@#1{\hbox{#1}}
\def\text@#1{\mathchoice
 {\hbox{\everymath{\displaystyle}\def\textfonti{\the\textfont1 }%
    \def\textfontii{\the\textfont2 }\textdef@@ T#1}}
 {\hbox{\everymath{\textstyle}\def\textfonti{\the\textfont1 }%
    \def\textfontii{\the\textfont2 }\textdef@@ T#1}}
 {\hbox{\everymath{\scriptstyle}\def\textfonti{\the\scriptfont1 }%
   \def\textfontii{\the\scriptfont2 }\textdef@@ S\rm#1}}
 {\hbox{\everymath{\scriptscriptstyle}\def\textfonti{\the\scriptscriptfont1 }%
   \def\textfontii{\the\scriptscriptfont2 }\textdef@@ s\rm#1}}}
\def\textdef@@#1{\textdef@#1\rm \textdef@#1\bf
   \textdef@#1\sl \textdef@#1\it}

\def\textdef@#1#2{\def\next{\csname\expandafter\eat@\string#2fam\endcsname}%
\if S#1\edef#2{\the\scriptfont\next\relax}%
 \else\if s#1\edef#2{\the\scriptscriptfont\next\relax}%
 \else\edef#2{\the\textfont\next\relax}\fi\fi}
\scriptfont\itfam=\tenit \scriptscriptfont\itfam=\tenit
\scriptfont\slfam=\tensl \scriptscriptfont\slfam=\tensl
\mathcode`\0="0030
\mathcode`\1="0031
\mathcode`\2="0032
\mathcode`\3="0033
\mathcode`\4="0034
\mathcode`\5="0035
\mathcode`\6="0036
\mathcode`\7="0037
\mathcode`\8="0038
\mathcode`\9="0039
\def\Cal{\relax\ifmmode\let\next=\Cal@\else
     \def\next{\errmessage{Use \string\Cal\space only in math mode}}\fi\next}
\def\Cal@#1{{\fam2 #1}}
\def\bold{\relax\ifmmode\let\next=\bold@\else
   \def\next{\errmessage{Use \string\bold\space only in math
      mode}}\fi\next}\def\bold@#1{{\fam\bffam #1}}
\mathchardef\Gamma="0000
\mathchardef\Delta="0001
\mathchardef\Theta="0002
\mathchardef\Lambda="0003
\mathchardef\Xi="0004
\mathchardef\Pi="0005
\mathchardef\Sigma="0006
\mathchardef\Upsilon="0007
\mathchardef\Phi="0008
\mathchardef\Psi="0009
\mathchardef\Omega="000A
\mathchardef\varGamma="0100
\mathchardef\varDelta="0101
\mathchardef\varTheta="0102
\mathchardef\varLambda="0103
\mathchardef\varXi="0104
\mathchardef\varPi="0105
\mathchardef\varSigma="0106
\mathchardef\varUpsilon="0107
\mathchardef\varPhi="0108
\mathchardef\varPsi="0109
\mathchardef\varOmega="010A
\font\dummyft@=dummy
\fontdimen1 \dummyft@=\z@
\fontdimen2 \dummyft@=\z@
\fontdimen3 \dummyft@=\z@
\fontdimen4 \dummyft@=\z@
\fontdimen5 \dummyft@=\z@
\fontdimen6 \dummyft@=\z@
\fontdimen7 \dummyft@=\z@
\fontdimen8 \dummyft@=\z@
\fontdimen9 \dummyft@=\z@
\fontdimen10 \dummyft@=\z@
\fontdimen11 \dummyft@=\z@
\fontdimen12 \dummyft@=\z@
\fontdimen13 \dummyft@=\z@
\fontdimen14 \dummyft@=\z@
\fontdimen15 \dummyft@=\z@
\fontdimen16 \dummyft@=\z@
\fontdimen17 \dummyft@=\z@
\fontdimen18 \dummyft@=\z@
\fontdimen19 \dummyft@=\z@
\fontdimen20 \dummyft@=\z@
\fontdimen21 \dummyft@=\z@
\fontdimen22 \dummyft@=\z@
\def\fontlist@{\\{\tenrm}\\{\sevenrm}\\{\fiverm}\\{\teni}\\{\seveni}%
 \\{\fivei}\\{\tensy}\\{\sevensy}\\{\fivesy}\\{\tenex}\\{\tenbf}\\{\sevenbf}%
 \\{\fivebf}\\{\tensl}\\{\tenit}\\{\tensmc}}
\def\dodummy@{{\def\\##1{\global\let##1=\dummyft@}\fontlist@}}
\newif\ifsyntax@
\newcount\countxviii@
\def\newtoks@{\alloc@5\toks\toksdef\@cclvi}
\def\nopages@{\output={\setbox\z@=\box\@cclv \deadcycles=\z@}\newtoks@\output}
\def\syntax{\syntax@true\dodummy@\countxviii@=\count18
\loop \ifnum\countxviii@ > \z@ \textfont\countxviii@=\dummyft@
   \scriptfont\countxviii@=\dummyft@ \scriptscriptfont\countxviii@=\dummyft@
     \advance\countxviii@ by-\@ne\repeat
\dummyft@\tracinglostchars=\z@
  \nopages@\frenchspacing\hbadness=\@M}
\def\magstep#1{\ifcase#1 1000\or
 1200\or 1440\or 1728\or 2074\or 2488\or
 \errmessage{\string\magstep\space only works up to 5}\fi\relax}
{\lccode`\2=`\p \lccode`\3=`\t
 \lowercase{\gdef\tru@#123{#1truept}}}

\def\scaletype#1{\mag=#1\relax
 \hsize=\expandafter\tru@\the\hsize
 \vsize=\expandafter\tru@\the\vsize
 \dimen\footins=\expandafter\tru@\the\dimen\footins}

\def\scalefont#1#2\andcallit#3{\edef\font@{\the\font}#1\font#3=
  \fontname\font\space scaled #2\relax\font@}
\def\Mag@#1#2{\ifdim#1<1pt\multiply#1 #2\relax\divide#1 1000 \else
  \ifdim#1<10pt\divide#1 10 \multiply#1 #2\relax\divide#1 100\else
  \divide#1 100 \multiply#1 #2\relax\divide#1 10 \fi\fi}
\def\scalelinespacing#1{\Mag@\baselineskip{#1}\Mag@\lineskip{#1}%
  \Mag@\lineskiplimit{#1}}
\def\wlog#1{\immediate\write-1{#1}}
\catcode`\@=\active

\input vanilla.sty
\magnification \magstep 1
\baselineskip = 12pt
\voffset = .2cm
\vglue 1.2cm
\hsize 13.34cm
\vsize 20.14cm
\parindent = 0mm
\font\big =cmb10 scaled \magstep 2
\font\smrm = cmr9

\TagsOnRight

\title \big
On certain organic compounds with one mono-substitution \\
\big and at least three di-substitution homogeneous derivatives
\footnote""{\smrm Research partially supported by Grant MM-1106/2001 of the
Bulgarian Foundation of Scientific Research}

\endtitle

\author
Valentin Vankov Iliev
\endauthor

\centerline {\it Section of Algebra, Institute of Mathematics and Informatics}
\centerline {\it Bulgarian Academy of Sciences, 1113 Sofia, Bulgaria}
\centerline {\it E-mail: viliev\@math.bas.bg, viliev\@aubg.bg}

\heading
1. Introduction
\endheading

A starting point of Lunn-Senior's theory of assigning a permutation group of
symmetry of degree $d$ to a given molecular structure divided into skeleton
and $d$ univalent substituents is the following old observation:  the number
of its substitution isomers does not depend on the nature of the ligants but
only on the numbers $\lambda_k$ of members of their different types $x_k$,
$k=1,2,\hdots$, and on the skeleton.  The only natural restriction is that if
the skeleton contains an univalent atom (or radical), then no univalent
substituent is to be identical with this atom (radical).  As far as the order
of ligants is irrelevant, we obtain a {\it partition}
$(\lambda_1,\lambda_2,\hdots,\lambda_d)$ of the number $d$, that is,
$\lambda_1\geq\lambda_2\geq\cdots\geq\lambda_d\geq 0$, and
$\lambda_1+\lambda_2+\cdots+\lambda_d=d$.  Plainly, the monomial
$$
x_1^{\lambda_1}x_2^{\lambda_2}\hdots x_d^{\lambda_d}
$$
is an exotic representation of substituents' {\it empirical formula} of the
molecular structure under question.  If $\Theta$ is the empirical formula of
the skeleton, then
$$
\Theta x_1^{\lambda_1}x_2^{\lambda_2}\hdots x_d^{\lambda_d}
$$
is the empirical formula of the molecule.  The additional information that
makes difference between its empirical and structural formulae consists of a
set of lists $A_k$, $k=1,2,\hdots, d$, each one enumerating the
unsatisfied valencies of the skeleton occupied by the identical ligants of type
$x_k$.  If a numeration $1,2,\hdots,d$ of the unsatisfied valencies is fixed
once and for all, then $A_k$ are simply pairwise disjoint subsets of the
integer-valued interval $[1,d]$, such that $[1,d]=\cup_kA_k$.  Thus, the
mathematical model of a {\it structural formula} of the substituents of a
molecular structure with empirical formula
$x_1^{\lambda_1}x_2^{\lambda_2}\hdots x_d^{\lambda_d}$, is a {\it tabloid
$A=(A_1,A_2,\hdots,A_d)$ with $d$ nodes of shape
$\lambda=(\lambda_1,\lambda_2,\hdots,\lambda_d)$}:
$$
A=
\left.
\matrix \format \c & \quad \c & \quad \c & \quad \c & \quad \c & \quad \c &
\quad \c & \quad \c \\
a_{1,1}, & a_{1,2}, & \hdots & \hdots & \hdots & a_{1,\lambda_1} &\hbox{\ \ } &
\hbox{\ the component \ } A_1  \cr
a_{2,1}, & a_{2,2}, & \hdots & \hdots & a_{2,\lambda_2} & \hbox{ \ \ } & \hbox{
\ \ } & \hbox{\ the component\ } A_2  \cr
\hbox{\ \ } & \hbox{\ \ } & \vdots & \hbox{\ \ } & \hbox{\ \ } & \hbox{\ \ }
& \hbox{\ \ } & \vdots   \cr
a_{t,1}, & \hdots & \hbox{\ \ } & \hbox{\ \ } & \hbox{\ \ } & \hbox{\ \ }  &
\hbox{\ \ } & \hbox{\ the component\ } A_t  \cr
\hbox{\ \ } & \hbox{\ \ } & \hbox{\ \ } & \hbox{\ \ } & \downarrow & \varphi &
\hbox{\ \ } & \hbox{\ \ }  \cr
\endmatrix
\right.
$$
$$
\lambda=
\left.
\matrix \format \c & \quad \c & \quad \c & \quad \c & \quad \c & \quad \c &
\quad \c & \quad \c \\
\times & \times & \cdots & \cdots & \cdots & \times &\hbox{\ \ } &
\lambda_1\hbox{\ nodes\ }   \cr
\times & \times & \cdots & \cdots & \times & \hbox{ \ \ } & \hbox{ \ \
} & \lambda_2\hbox{\ nodes\ }   \cr
\hbox{\ \ } & \hbox{\ \ } & \vdots & \hbox{\ \ } & \hbox{\ \ } & \hbox{\ \ }
& \hbox{\ \ } & \vdots   \cr
\times & \cdots & \hbox{\ \ } & \hbox{\ \ } & \hbox{\ \ } & \hbox{\ \ }  &
\hbox{\ \ } & \lambda_t\hbox{\ nodes\ }   \cr
\endmatrix
\right.
$$
Here $\varphi\colon T_d\to P_d$ is the natural projection of the set
$T_d$ of all tabloids with $d$ nodes onto the set $P_d$ of all partitions of
$d$ that maps the tabloid $A$ onto its {\it shape} $\lambda$:
$\lambda_1=|A_1|$, $\lambda_2=|A_2|,\hdots$, $\lambda_d=|A_d|$.

The structural formula of a molecule encodes its ``connexity data", and does
not reflect in full so called ``space configuration", because the latter is a
special representation of the former.  ``Connexity is a relation of order
independent of considerations of space.  The ``structural" relations treated by
chemists are relations of just this sort, and it is unfortunate that the word
structure as used by engineers, etc., should carry with it geometrical
connotations which are too special for chemistry" [6, p.  1030].

The inverse image $T_\lambda=\varphi^{-1}(\lambda)$ consists of all structural
formulae of the substituents with empirical formula
$x_1^{\lambda_1}x_2^{\lambda_2}\hdots x_d^{\lambda_d}$. The fibers $T_\lambda$,
$\lambda\in P_d$, of the map $\varphi$ are the stages where the drama of
isomerism is performed.

In [6], Lunn and Senior build in the phenomenon of isomerism of a certain type
in the above mathematical model by means of action of a symmetry group $G$,
consisting of permutations of the $d$ unsatisfied valencies of the skeleton, and
such that any isomer of the given empirical formula $\Theta
x_1^{\lambda_1}x_2^{\lambda_2}\hdots x_d^{\lambda_d}$ is represented by a
$G$-orbit in $T_\lambda$. The group $G$ acts on the set $T_d$ of structural
formulae by the rule
$$
\sigma(A_1,A_2,\hdots,A_d)=
(\sigma(A_1),\sigma(A_2),\hdots,\sigma(A_d)),
$$
and produces the spaces $T_{\lambda;G}$ of $G$-orbits of the structural
formulae from $T_\lambda$.  The number $n_{\lambda;G}$ of these $G$-orbits is
therefore an upper bound for the number $N_{\lambda;\Theta}$ of experimentally
known derivatives with composition $\Theta x_1^{\lambda_1}x_2^{\lambda_2}\hdots
x_d^{\lambda_d}$:
$$
N_{\lambda;\Theta}\leq n_{\lambda;G}
$$
for any partition $\lambda\in P_d$.  In the cases of mono-substituted
derivatives ($\lambda=(d-1,1)$), di-substituted homogeneous derivatives
($\lambda=(d-2,2)$), and di-substituted heterogeneous derivatives
($\lambda=(d-2,1^2)$), the experimenters, sometimes, are certain that the
corresponding numbers $N_{\lambda;\Theta}$ attain their maximum values
$n_{\lambda;G}$.  In other words, all possible {\it $\lambda$-derivatives} are
prepared.  In the ideal (but unattainable) situation
$N_{\lambda;\Theta}=n_{\lambda;G}$ for all partitions $\lambda\in P_d$, and
these equalities define the symmetry group $G$ up to so called combinatorial
equivalence (See [6, IV], [7, 26], [2, 5.2.5]).

The {\it simple substitution reactions}
$$
x_1^{\mu_1}\hdots x_i^{\mu_i}\hdots x_j^{\mu_j}\hdots\longrightarrow
x_1^{\lambda_1}\hdots x_i^{\lambda_i}\hdots x_j^{\lambda_j}\hdots,
$$
where $\lambda, \mu\in P_d$, and $\mu_1=\lambda_1,\hdots$,
$\mu_i=\lambda_i+1,\hdots$, $\mu_j=\lambda_j-1,\hdots$, $\mu_d=\lambda_d$, that
is, the replacement of a ligant of type $x_i$ by a ligant of type $x_j$, $j<i$,
are encoded in the mathematical model via two partial orderings: on the level
of empirical formulae we write $\lambda <\mu$,  and on the level of the
structural picture
$$
B=(B_1,B_2,\hdots,B_d)\longrightarrow
A=(A_1,A_2,\hdots,A_d)
$$
$A, B\in T_d$, $\lambda=\varphi(A)$, $\mu=\varphi(B)$, of the above simple
substitution reaction, where $A$ is obtained from $B$ by moving an element
$s\in B_i$ in the set $B_j$, we write $A<B$.  More generally, we write
$\lambda <\mu$ if $\lambda$ can be got from $\mu$ by a finite number of the
above simple substitutions (this is the well known {\it dominance order} of
partitions, See [5, 6.1]), and we write $A<B$ if $A$ can be obtained from $B$
via a finite sequence of the above simple movements of elements (See [2, 3.2].
The latter ordering can be pulled down on the orbit-space $T_{d;G}=G\backslash
T_d$:  $a<b$ if there are $A\in a$, $B\in b$ with $A<B$ (See [2, 4.1]).  If
$a<b$, $a,b\in T_{d;G}$, the product which corresponds to $a$ can, in
principle, be synthesized from the product which corresponds to $b$ via a
finite sequence of simple substitution reactions.  Thus, the partially ordered
set $T_{d;G}$ portrays the possible {\it genetic relations} among the
derivatives of the molecule under consideration (See [2]).

In this paper we consider parent substances with molecules that can be divided
into a skeleton and six univalent substituents, and have the properties
mentioned in the title.  Two instances are the molecules of benzene $C_6H_6$
and cyclopropane $C_3H_6$, which have one mono-substitution derivative, and
three and four di-substitution homogeneous derivatives, respectively.

The paper is stratified as follows. In Section 2, Theorem 2.1 describes the
Lunn-Senior's group $G$ of substitution isomerism of our compounds and
Corollaries 2.7 and 2.8 give upper bounds of the numbers of their
di-substitution, and tri-substitution homogeneous derivatives.  In Sections 3,
4, and 5, we list the possible simple substitution reactions among
di-substitution homogeneous derivatives, on one hand, and di-substitution
heterogeneous, and tri-substitution homogeneous derivatives, on the other.
These substitution reactions allow us to identify some derivatives with their
structural formulae.

\heading
2. The Lunn-Senior's Group of Substitution Isomerism

\endheading

The theorem below gives a characterization of the Lunn-Senior's groups of
substitution isomerism of the compounds from the title.

\proclaim{Theorem 2.1} If an organic compound consists of a skeleton with six
univalent substituents and has one mono-substitution and at least three
di-substitution homogeneous derivatives, then its Lunn-Senior's group of
substitution isomerism is conjugated in $S_6$ either to the dihedral group
$$
\langle (123)(456), (14)(26)(35), (14)(25)(36)\rangle
$$
of order $12$, or to the cyclic group
$$
\langle (123456)\rangle
$$
of order $6$, or to the dihedral group
$$
\langle (123)(456), (14)(26)(35)\rangle
$$
of order $6$.

\endproclaim

\demo{Proof} Since there exists only one mono-substitution derivative, we have
$n_{\left(5,1\right);G}=1$, so the Lunn-Senior's group $G\leq S_6$ of
substitution isomerism is transitive (See [3, 3.1.1].  The existence of at
least three di-substitution derivatives means that
$$
n_{\left(4,2\right);G}\geq 3.  \tag 2.2
$$
Since the partition $(4,2)$ dominates the partition $(4,1^2)$ with respect to
the dominance order,
[2, (5.3.2)] implies
$$
n_{\left(4,2\right);G}\leq n_{\left(4,1^2\right);G}. \tag 2.3
$$
In particular, $n_{\left(4,1^2\right);G}\geq 3$. Therefore [4, (6.1.1)]
and [4,(6.1.2)] yield $g_{\left(4,2\right);G}= g_{\left(4,1^2\right);G}=
g_{\left(3,2,1\right);G}= g_{\left(2,1^4\right);G}=g_{\left(3,1^3\right);G}=0$.
Then the linear system [3, (3.2.1)] becomes
$$
\left.
\matrix
\format \c & \ \c & \ \c &\ \c &\ \c &\ \c & \ \c & \ \c & \ \c & \ \c &\ \c &\
\c &\ \c &\ \c & \ \c &\ \c &\ \c & \c \\
g_{\left(6\right);G} & + & g_{\left(3^2\right);G} & + & g_{\left(2^3\right);G} &
+ & g_{\left(2^2,1^2\right);G} & - & (|G|-1) & = & 0 \cr
\hbox{ } & \hbox{ } & 2g_{\left(3^2\right);G} & \hbox{ }  & \hbox{ }  &
+ & 4g_{\left(2^2,1^2\right);G} & - & (|G|n_{\left(3^2\right);G}-20) & = & 0 \cr
\hbox{ } & \hbox{ } & \hbox{ } & \hbox{ }  & 3g_{\left(2^3\right);G} &
+ & 3g_{\left(2^2,1^2\right);G} & - & (|G|n_{\left(4,2\right);G}-15) & = & 0 \cr
\hbox{ } & \hbox{ } & \hbox{ } & \hbox{ }  & 6g_{\left(2^3\right);G} &
+ & 6g_{\left(2^2,1^2\right);G} & - & (|G|n_{\left(2^3\right);G}-90) & = & 0 \cr
\hbox{ } & \hbox{ } & \hbox{ } & \hbox{ }  & \hbox{ } & \hbox{ }
& 2g_{\left(2^2,1^2\right);G} & - & (|G|n_{\left(5,1\right);G}-6) & = & 0 \cr
\hbox{ } & \hbox{ } & \hbox{ } & \hbox{ }  & \hbox{ } & \hbox{ }
& 2g_{\left(2^2,1^2\right);G} & - & (|G|n_{\left(4,1^2\right);G}-30) & = & 0 \cr
\hbox{ } & \hbox{ } & \hbox{ } & \hbox{ }  & \hbox{ } & \hbox{ }
& 4g_{\left(2^2,1^2\right);G} & - & (|G|n_{\left(3,2,1\right);G}-60) & = & 0 \cr
\hbox{ } & \hbox{ } & \hbox{ } & \hbox{ }  & \hbox{ } & \hbox{ }
& 4g_{\left(2^2,1^2\right);G} & - & (|G|n_{\left(2^2,1^2\right);G}-180) & = & 0
\cr
\hbox{ } & \hbox{ } & \hbox{ } & \hbox{ } & \hbox{ } & \hbox{ } & \hbox{ } &
\hbox{ } & (|G|n_{\left(3,1^3\right);G}-120) & = & 0 \cr
\hbox{ } & \hbox{ } & \hbox{ } & \hbox{ } & \hbox{ } & \hbox{ } & \hbox{ } &
\hbox{ } & (|G|n_{\left(2,1^4\right);G}-360) & = & 0 \cr
\hbox{ } & \hbox{ } & \hbox{ } & \hbox{ } & \hbox{ } & \hbox{ } & \hbox{ } &
\hbox{ } & (|G|n_{\left(1^6\right);G}-720) & = & 0 \cr
\endmatrix
\right.\tag 2.4
$$

Since $n_{\left(5,1\right);G}=1$, the fifth and sixth equalities of the system
(2.4) yield $|G|\geq 6$, and
$$
|G|(n_{\left(4,1^2\right);G}-1)=24. \tag 2.5
$$
Then the inequality $n_{\left(4,1^2\right);G}\geq 3$ implies $|G|\leq 12$.
On the other hand, the third and fourth equalities of (2.4) imply
$$
|G|(n_{\left(2^3\right);G}-2n_{\left(4,2\right);G})=60.
$$
Thus, $|G|$ is a common divisor of $24$ and $60$, so we obtain two
possibilities for the order $|G|$ of the group $G$:  $|G|=12$ or $|G|=6$.

If $|G|=12$, then from (2.5) we get $n_{\left(4,1^2\right);G}=3$, and the
inequalities (2.2) and (2.3) yield $n_{\left(4,2\right);G}=3$.  Now, the
fifth and the third equality of (2.4) imply $g_{\left(2^2,1^2\right);G}=3$
and $g_{\left(2^3\right);G}=4$.  Hence the first equality of (2.4) yields
$g_{\left(6\right);G} + g_{\left(3^2\right);G}=4$.  The equality
$g_{\left(3^2\right);G}=0$ is impossible since for every cycle $\sigma\in G$ of
length $6$ its square $\sigma^2$ has cyclic type $(3^2)$.  Therefore
$g_{\left(6\right);G}=g_{\left(3^2\right);G}=2$.  Let $\sigma$ be a cycle
of length $6$. After eventual conjugation,
we can suppose that $\sigma^2=(123)(456)\in G$.  Now, consider the cyclic group
$K=\langle\sigma\rangle$ of order $6$ and its cyclic subgroup $H=\langle
(123)(456)\rangle$ that contains the two elements of $G$ of cyclic type
$(3^2)$.  If $\iota$ is one of the the elements of $G$ of cyclic type $(2^3)$,
then
$$
\iota H\iota^{-1}=H, \tag 2.6
$$
so $L=H\langle\iota\rangle$ is a subgroup of $G$ of order $6$.  Now, we choose
$\iota\notin K$ (since $g_{\left(2^3\right);G}=4$ there are three elements of
cyclic type $(2^3)$ outside $K$).  If we suppose that $L$ is cyclic, then we
would have $L=K$ (the two elements of order $6$ in $G$ are in $K$), and in
particular, $\iota\in K$:  a contradiction.  Hence $L$ is isomorphic to the
dihedral group of order $6$.  Further, the equality (2.6) and the
considerations in [4, 7.1] yield that we can set $\iota=(14)(26)(35)$, so
$L=\langle (123)(456), (14)(26)(35)\rangle.$ Now, in accord to [4, 7.3.1],
we get that the group $G$ is conjugated to the group $\langle (123)(456),
(14)(26)(35), (14)(25)(36)\rangle.$

If $|G|=6$, then $n_{\left(4,1^2\right);G}=5$ and the fifth equality of (2.4)
implies $g_{\left(2^2,1^2\right);G}=0$.  Then the first equality of the system
(2.4) becomes $g_{\left(6\right);G} + g_{\left(3^2\right);G} +
g_{\left(2^3\right);G}=5.$

If $G$ is the cyclic group of order $6$, then it is generated, up to
conjugation, by the cycle $(123456)$, and $g_{\left(6\right);G}=2$,
$g_{\left(3^2\right);G}=2$, and $g_{\left(2^3\right);G}=1$.  Now, the third
equality of (2.4) yields $n_{\left(4,2\right);G}=3$.

If $G$ is the dihedral group of order $6$, then $g_{\left(3^2\right);G}=2$,
$g_{\left(2^3\right);G}=3$, and in accordance to the third equality of (2.4),
we obtain $n_{\left(4,2\right);G}=4$.  Now, we apply [4, 5.1.1].

\endproclaim

Theorem 2.1 implies immediately the following two corollaries which yield the
numbers of derivatives of the molecules under consideration.

\proclaim{Corollary 2.7} If an organic compound consists of a skeleton with
six univalent substituents and has one mono-substitution and at least three
di-substitution homogeneous derivatives, and if its Lunn-Senior's group of
substitution isomerism has order $12$, then this compound has exactly three
di-substitution homogeneous derivatives, at most three di-substitution
heterogeneous derivatives, and at most three tri-substitution homogeneous
derivatives.

\endproclaim

\proclaim{Corollary 2.8} If an organic compound consists of a skeleton with
six univalent substituents and has one mono-substitution and at least three
di-substitution homogeneous derivatives, and if its Lunn-Senior's
group $G$ of substitution isomerism has order $6$, then this compound has
exactly three di-substitution homogeneous derivatives, at most five
di-substitution heterogeneous derivatives, and at most four tri-substitution
homogeneous derivatives in case $G$ is cyclic, or has three or four
di-substitution homogeneous derivatives, at most five di-substitution
heterogeneous derivatives, and at most four tri-substitution homogeneous
derivatives in case $G$ is dihedral.

\endproclaim

\heading
3. Genetic Relations: the Group $G$ has Order $12$

\endheading

Here we consider the possible genetic relations among the derivatives of our
molecule structure in the case when its Lunn-Senior's group $G$ of substitution
isomerism has order $12$.  An example is the benzen molecule $C_6H_6$ (See [6,
VI] or [1], or [2, 6.3]).  In accord to [2, 6.3] and Theorem 3.1, we may
suppose $G=\langle (123456), (13)(46)\rangle$ and then we obtain
$T_{\left(4,2\right);G}=\{a_{\left(4,2\right)},b_{\left(4,2\right)},
c_{\left(4,2\right)}\}$, where:

$a_{\left(4,2\right)}$ is the $G$-orbit
$$
\{(\{1,2,4,5\},\{3,6\}),
(\{2,3,5,6\},\{1,4\}),
(\{1,3,4,6\},\{2,5\})\}
$$
of the tabloid $A^{\left(4,2\right)}=(\{1,2,4,5\},\{3,6\})$,

$b_{\left(4,2\right)}$ is the $G$-orbit
$$
\{(\{1,2,3,4\},\{5,6\}),
(\{2,3,4,5\},\{1,6\}),
(\{3,4,5,6\},\{1,2\}),
$$
$$
(\{1,4,5,6\},\{2,3\}),
(\{1,2,5,6\},\{3,4\}),
(\{1,2,3,6\},\{4,5\})\}
$$
of the tabloid
$B^{\left(4,2\right)}=
(\{1,2,3,4\},\{5,6\})$,

$c_{\left(4,2\right)}$ is the $G$-orbit
$$
\{(\{1,2,4,6\},\{3,5\}),
(\{1,2,3,5\},\{4,6\}),
(\{2,3,4,6\},\{1,5\}),
$$
$$
(\{1,3,4,5\},\{2,6\}),
(\{2,4,5,6\},\{1,3\}),
(\{1,3,5,6\},\{2,4\})\}
$$
of the tabloid
$C^{\left(4,2\right)}=
(\{1,2,4,6\},\{3,5\})$.

Further, we get
$T_{\left(3^2\right);G}=\{a_{\left(3^2\right)},b_{\left(3^2\right)},
c_{\left(3^2\right)}\}$, where:

$a_{\left(3^2\right)}$ is the $G$-orbit
$$
\{(\{1,2,4\},\{3,5,6\}),
(\{2,3,5\},\{1,4,6\}),
(\{3,4,6\},\{1,2,5\}),
(\{1,4,5\},\{2,3,6\}),
$$
$$
(\{2,5,6\},\{1,3,4\}),
(\{1,3,6\},\{2,4,5\}),
(\{2,3,6\},\{1,4,5\}),
(\{1,2,5\},\{3,4,6\}),
$$
$$
(\{1,4,6\},\{2,3,5\}),
(\{3,5,6\},\{1,2,4\}),
(\{2,4,5\},\{1,3,6\}),
(\{1,3,4\},\{2,5,6\})\}
$$
of the tabloid
$A^{\left(3^2\right)}=
(\{1,2,4\},\{3,5,6\})$;

$b_{\left(3^2\right)}$ is the $G$-orbit
$$
\{(\{1,2,3\},\{4,5,6\}),
(\{2,3,4\},\{1,5,6\}),
(\{3,4,5\},\{1,2,6\}),
$$
$$
(\{4,5,6\},\{1,2,3\}),
(\{1,5,6\},\{2,3,4\}),
(\{1,2,6\},\{3,4,5\})\}
$$
of the tabloid
$B^{\left(3^2\right)}=
(\{1,2,3\},\{4,5,6\})$;

$c_{\left(3^2\right)}$ is the $G$-orbit
$$
\{(\{1,3,5\},\{2,4,6\}),
(\{2,4,6\},\{1,3,5\})\}
$$
of the tabloid
$C^{\left(3^2\right)}=
(\{1,3,5\},\{2,4,6\})$.

Moreover, we obtain

$T_{\left(4,1^2\right);G}=\{a_{\left(4,1^2\right)},b_{\left(4,1^2\right)},
c_{\left(4,1^2\right)}\}$, where:

$a_{\left(4,1^2\right)}$ is the $G$-orbit
$$
\{(\{1,2,4,5\},\{3\},\{6\}),
(\{2,3,5,6\},\{4\},\{1\}),
(\{1,3,4,6\},\{5\},\{2\}),
$$
$$
(\{1,2,4,5\},\{6\},\{3\}),
(\{2,3,5,6\},\{1\},\{4\}),
(\{1,3,4,6\},\{2\},\{5\})\}
$$
of the tabloid $A^{\left(4,1^2\right)}=(\{1,2,4,5\},\{3\},\{6\})$,

$b_{\left(4,1^2\right)}$ is the $G$-orbit
$$
\{(\{1,2,3,4\},\{5\},\{6\}),
(\{2,3,4,5\},\{6\},\{1\}),
(\{3,4,5,6\},\{1\},\{2\}),
$$
$$
(\{1,4,5,6\},\{2\},\{3\}),
(\{1,2,5,6\},\{3\},\{4\}),
(\{1,2,3,6\},\{4\},\{5\}),
$$
$$
(\{1,2,3,4\},\{6\},\{5\}),
(\{2,3,4,5\},\{1\},\{6\}),
(\{3,4,5,6\},\{2\},\{1\}),
$$
$$
(\{1,4,5,6\},\{3\},\{2\}),
(\{1,2,5,6\},\{4\},\{3\}),
(\{1,2,3,6\},\{5\},\{4\})\}
$$

of the tabloid
$B^{\left(4,1^2\right)}=
(\{1,2,3,4\},\{5\},\{6\})$,

$c_{\left(4,1^2\right)}$ is the $G$-orbit
$$
\{(\{1,2,4,6\},\{3\},\{5\}),
(\{1,2,3,5\},\{4\},\{6\}),
(\{2,3,4,6\},\{5\},\{1\}),
$$
$$
(\{1,3,4,5\},\{6\},\{2\}),
(\{2,4,5,6\},\{3\},\{1\}),
(\{1,3,5,6\},\{4\},\{2\}),
$$
$$
(\{1,2,4,6\},\{5\},\{3\}),
(\{1,2,3,5\},\{6\},\{4\}),
(\{2,3,4,6\},\{1\},\{5\}),
$$
$$
(\{1,3,4,5\},\{2\},\{6\}),
(\{2,4,5,6\},\{1\},\{3\}),
(\{1,3,5,6\},\{2\},\{4\})\}
$$
of the tabloid
$C^{\left(4,1^2\right)}=
(\{1,2,4,6\},\{3\},\{5\})$.

Since
$$
A^{\left(3^2\right)}<A^{\left(4,2\right)},\hbox{\ }
A^{\left(3^2\right)}<B^{\left(4,2\right)},\hbox{\ }
A^{\left(3^2\right)}<C^{\left(4,2\right)},
$$
$$
B^{\left(3^2\right)}<B^{\left(4,2\right)},\hbox{\ }
B^{\left(3^2\right)}<(123456)C^{\left(4,2\right)},\hbox{\ }
C^{\left(3^2\right)}<(123456)C^{\left(4,2\right)},
$$
and since
$$
A^{\left(4,1^2\right)}<A^{\left(4,2\right)},\hbox{\ }
B^{\left(4,1^2\right)}<B^{\left(4,2\right)},\hbox{\ }
C^{\left(4,1^2\right)}<C^{\left(4,2\right)},
$$
we have the following inequalities
$$
a_{\left(3^2\right)} < a_{\left(4,2\right)},\hbox{\ }
a_{\left(3^2\right)} < b_{\left(4,2\right)},\hbox{\ }
a_{\left(3^2\right)} < c_{\left(4,2\right)},
$$
$$
b_{\left(3^2\right)} < b_{\left(4,2\right)},\hbox{\ }
b_{\left(3^2\right)} < c_{\left(4,2\right)},\hbox{\ }
c_{\left(3^2\right)} < c_{\left(4,2\right)},
$$
and
$$
a_{\left(4,1^2\right)} < a_{\left(4,2\right)},\hbox{\ }
a_{\left(4,1^2\right)} < b_{\left(4,2\right)},\hbox{\ }
a_{\left(4,1^2\right)} < c_{\left(4,2\right)}.
$$

The diagrams below represent \lq\lq K\"orner like" relations between
homogeneous di- and tri-substitution products of our molecule structure, which
can be used for complete identification of these six derivatives.
$$
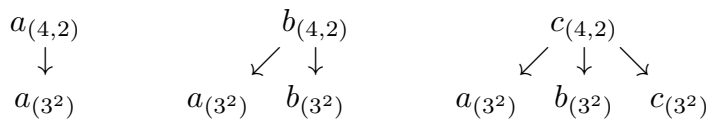
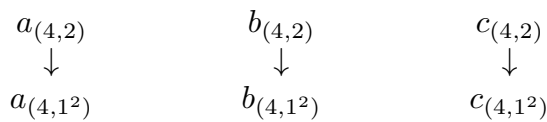
\left.
\matrix \format \c & \c & \c & \c & \c & \c & \c & \c & \c & \c & \c & \c & \c
& \c &  \c  & \c & \c & \c & \c \\
a_{\left(4,2\right)}  & \hbox{\ \ } & \hbox{\ \ } & \hbox{\ \ } & \hbox{\ \ }
& \hbox{\ \ } & \hbox{\ \ } & \hbox{\ \ } & b_{\left(4,2\right)}  & \hbox{\ \
} & \hbox{\ \ } & \hbox{\ \ } & \hbox{\ \ } & \hbox{\ \ } & \hbox{\ \ } &
\hbox{\ \ } & c_{\left(4,2\right)} & \hbox{\ \ } & \hbox{\ \ }  \cr
\downarrow & \hbox{\ \ } & \hbox{\ \ } & \hbox{\ \ } & \hbox{\ \ } & \hbox{\ \
} & \hbox{\ \ } & \swarrow & \downarrow & \hbox{\ \ } & \hbox{\ \ } &
\hbox{\ \ } & \hbox{\ \ } & \hbox{\ \ } & \hbox{\ \ } & \swarrow & \downarrow &
\searrow & \hbox{\ \ } \cr
a_{\left(3^2\right)} & \hbox{\ \ } &\hbox{\ \ } & \hbox{\ \ } & \hbox{\ \ } &
\hbox{\ \ } & a_{\left(3^2\right)} & \hbox{\ \ } & b_{\left(3^2\right)} &
\hbox{\ \ } & \hbox{\ \ } & \hbox{\ \ } & \hbox{\ \ } & \hbox{\ \ } &
a_{\left(3^2\right)} & \hbox{\ \ } & b_{\left(3^2\right)} & \hbox{\ \ } &
c_{\left(3^2\right)} \cr
\endmatrix \right.
$$
The diagrams
$$
\left.
\matrix \format \c & \c & \c & \c & \c & \c & \c & \c & \c & \c & \c & \c & \c
& \c &  \c  & \c & \c & \c & \c \\
a_{\left(4,2\right)}  & \hbox{\ \ } & \hbox{\ \ } & \hbox{\ \ } & \hbox{\ \ }
& \hbox{\ \ } & \hbox{\ \ } & \hbox{\ \ } & b_{\left(4,2\right)}  & \hbox{\ \
} & \hbox{\ \ } & \hbox{\ \ } & \hbox{\ \ } & \hbox{\ \ } & \hbox{\ \ } &
\hbox{\ \ } & c_{\left(4,2\right)} & \hbox{\ \ } & \hbox{\ \ }  \cr
\downarrow & \hbox{\ \ } & \hbox{\ \ } & \hbox{\ \ } & \hbox{\ \ } & \hbox{\ \
} & \hbox{\ \ } & \hbox{\ \ } & \downarrow & \hbox{\ \ } & \hbox{\ \ } &
\hbox{\ \ } & \hbox{\ \ } & \hbox{\ \ } & \hbox{\ \ } & \hbox{\ \ } &
\downarrow & \hbox{\ \ } & \hbox{\ \ } \cr
a_{\left(4,1^2\right)} & \hbox{\ \ } &\hbox{\ \ } & \hbox{\ \ } & \hbox{\ \ } &
\hbox{\ \ } & \hbox{\ \ } & \hbox{\ \ } & b_{\left(4,1^2\right)} &
\hbox{\ \ } & \hbox{\ \ } & \hbox{\ \ } & \hbox{\ \ } & \hbox{\ \ } &
\hbox{\ \ } & \hbox{\ \ } & c_{\left(4,1^2\right)} & \hbox{\ \ } &
\hbox{\ \ } \cr
\endmatrix \right.
$$
show that, as a consequence, the heterogeneous di-substitution derivatives can
also be identified completely.

Here the arrow $a\rightarrow b$ means that $a>b$ and the product that
corresponds to $b$ can be obtained from the product that corresponds to $a$ via
a simple substitution reaction.

\heading
4. Genetic Relations: The Group $G$ has Order $6$ and is Cyclic

\endheading

In this section we describe the genetic relations of the molecule structure
under question when its Lunn-Senior's group $G$ of substitution isomerism is
cyclic of order $6$.  In accord with Theorem 2.1, we can suppose $G=\langle
(123456)\rangle$.  Then
$T_{\left(4,2\right);G}=\{a_{\left(4,2\right)},b_{\left(4,2\right)},
c_{\left(4,2\right)}\}$, where:

$a_{\left(4,2\right)}$ is the $G$-orbit
$$
\{(\{1,2,4,5\},\{3,6\}),
(\{2,3,5,6\},\{1,4\}),
(\{1,3,4,6\},\{2,5\})\},
$$
of the tabloid $A^{\left(4,2\right)}=(\{1,2,4,5\},\{3,6\})$,

$b_{\left(4,2\right)}$ is the $G$-orbit
$$
\{(\{1,2,3,4\},\{5,6\}),
(\{2,3,4,5\},\{1,6\}),
(\{3,4,5,6\},\{1,2\}),
$$
$$
(\{1,4,5,6\},\{2,3\}),
(\{1,2,5,6\},\{3,4\}),
(\{1,2,3,6\},\{4,5\})\}
$$
of the tabloid
$B^{\left(4,2\right)}=
(\{1,2,3,4\},\{5,6\})$,

$c_{\left(4,2\right)}$ is the $G$-orbit
$$
\{(\{1,2,4,6\},\{3,5\}),
(\{1,2,3,5\},\{4,6\}),
(\{2,3,4,6\},\{1,5\}),
$$
$$
(\{1,3,4,5\},\{2,6\}),
(\{2,4,5,6\},\{1,3\}),
(\{1,3,5,6\},\{2,4\})\},
$$
of the tabloid
$C^{\left(4,2\right)}=(\{1,2,4,6\},\{3,5\})$.

We have
$T_{\left(3^2\right);G}=\{a_{\left(3^2\right)},b_{\left(3^2\right)},
c_{\left(3^2\right)},d_{\left(3^2\right)}\}$, where:

$a_{\left(3^2\right)}$ is the $G$-orbit
$$
\{(\{1,2,4\},\{3,5,6\}),
(\{2,3,5\},\{1,4,6\}),
(\{3,4,6\},\{1,2,5\}),
$$
$$
(\{1,4,5\},\{2,3,6\}),
(\{2,5,6\},\{1,3,4\}),
(\{1,3,6\},\{2,4,5\})\},
$$
of the tabloid
$A^{\left(3^2\right)}=
(\{1,2,4\},\{3,5,6\})$;

$b_{\left(3^2\right)}$ is the $G$-orbit
$$
\{(\{1,2,5\},\{3,4,6\})\},
(\{2,3,6\},\{1,4,5\}),
(\{1,3,4\},\{2,5,6\}),
$$
$$
(\{2,4,5\},\{1,3,6\}),
(\{3,5,6\},\{1,2,4\}),
(\{1,4,6\},\{2,3,5\})\},
$$
of the tabloid
$B^{\left(3^2\right)}=
(\{1,2,5\},\{3,4,6\})$;

$c_{\left(3^2\right)}$ is the $G$-orbit
$$
\{(\{1,2,3\},\{4,5,6\}),
(\{2,3,4\},\{1,5,6\}),
(\{3,4,5\},\{1,2,6\}),
$$
$$
(\{4,5,6\},\{1,2,3\}),
(\{1,5,6\},\{2,3,4\}),
(\{1,2,6\},\{3,4,5\})\}
$$
of the tabloid
$C^{\left(3^2\right)}=
(\{1,2,3\},\{4,5,6\})$;

$d_{\left(3^2\right)}$ is the $G$-orbit

$$
\{(\{1,3,5\},\{2,4,6\}),
(\{2,4,6\},\{1,3,5\})\}
$$
of the tabloid
$D^{\left(3^2\right)}=
(\{1,3,5\},\{2,4,6\})$.

We also obtain
$T_{\left(4,1^2\right);G}=\{a_{\left(4,1^2\right)},b_{\left(4,1^2\right)},
c_{\left(4,1^2\right)},d_{\left(4,1^2\right)},e_{\left(4,1^2\right)}\}$, where:

$a_{\left(4,1^2\right)}$ is the $G$-orbit
$$
\{(\{1,2,4,5\},\{3\},\{6\}),
(\{2,3,5,6\},\{4\},\{1\}),
(\{1,3,4,6\},\{5\},\{2\}),
$$
$$
(\{1,2,4,5\},\{6\},\{3\}),
(\{2,3,5,6\},\{1\},\{4\}),
(\{1,3,4,6\},\{2\},\{5\})\}
$$
of the tabloid $A^{\left(4,1^2\right)}=(\{1,2,4,5\},\{3\},\{6\})$;

$b_{\left(4,1^2\right)}$ is the $G$-orbit
$$
\{(\{1,2,3,4\},\{5\},\{6\}),
(\{2,3,4,5\},\{6\},\{1\}),
(\{3,4,5,6\},\{1\},\{2\}),
$$
$$
(\{1,4,5,6\},\{2\},\{3\}),
(\{1,2,5,6\},\{3\},\{4\}),
(\{1,2,3,6\},\{4\},\{5\})\}
$$
of the tabloid
$B^{\left(4,1^2\right)}=
(\{1,2,3,4\},\{5\},\{6\})$;

$c_{\left(4,1^2\right)}$ is the $G$-orbit
$$
\{(\{1,2,3,4\},\{6\},\{5\}),
(\{2,3,4,5\},\{1\},\{6\}),
(\{3,4,5,6\},\{2\},\{1\}),
$$
$$
(\{1,4,5,6\},\{3\},\{2\}),
(\{1,2,5,6\},\{4\},\{3\}),
(\{1,2,3,6\},\{5\},\{4\})\}
$$
of the tabloid
$C^{\left(4,1^2\right)}=
(\{1,2,3,4\},\{6\},\{5\})$;

$d_{\left(4,1^2\right)}$ is the $G$-orbit
$$
\{(\{1,2,4,6\},\{3\},\{5\}),
(\{1,2,3,5\},\{4\},\{6\}),
(\{2,3,4,6\},\{5\},\{1\}),
$$
$$
(\{1,3,4,5\},\{6\},\{2\}),
(\{2,4,5,6\},\{1\},\{3\}),
(\{1,3,5,6\},\{2\},\{4\})\}
$$
of the tabloid
$D^{\left(4,1^2\right)}=(\{1,2,4,6\},\{3\},\{5\})$;

$e_{\left(4,1^2\right)}$ is the $G$-orbit
$$
\{(\{1,2,4,6\},\{5\},\{3\}),
(\{1,2,3,5\},\{6\},\{4\}),
(\{2,3,4,6\},\{1\},\{5\}),
$$
$$
(\{1,3,4,5\},\{2\},\{6\}),
(\{2,4,5,6\},\{3\},\{1\}),
(\{1,3,5,6\},\{4\},\{2\})\}
$$
of the tabloid
$E^{\left(4,1^2\right)}=(\{1,2,4,6\},\{5\},\{3\})$.

We have
$$
A^{\left(3^2\right)}<A^{\left(4,2\right)}, \hbox{\ }
A^{\left(3^2\right)}<B^{\left(4,2\right)},\hbox{\ }
A^{\left(3^2\right)}<C^{\left(4,2\right)},\hbox{\ }
$$
$$
B^{\left(3^2\right)}<A^{\left(4,2\right)},\hbox{\ }
B^{\left(3^2\right)}<(153)(264)B^{\left(4,2\right)},\hbox{\ }
B^{\left(3^2\right)}<(123456)C^{\left(4,2\right)},\hbox{\ }
$$
$$
C^{\left(3^2\right)}<B^{\left(4,2\right)},\hbox{\ }
C^{\left(3^2\right)}<(123456)C^{\left(4,2\right)},\hbox{\ }
D^{\left(3^2\right)}<(123456)C^{\left(4,2\right)},
$$
and
$$
A^{\left(4,1^2\right)}<A^{\left(4,2\right)},\hbox{\ }
B^{\left(4,1^2\right)}<B^{\left(4,2\right)},\hbox{\ }
C^{\left(4,1^2\right)}<B^{\left(4,2\right)},
$$
$$
D^{\left(4,1^2\right)}<C^{\left(4,2\right)},\hbox{\ }
E^{\left(4,1^2\right)}<C^{\left(4,2\right)},
$$
so
$$
a_{\left(3^2\right)}<a_{\left(4,2\right)}, \hbox{\ }
a_{\left(3^2\right)}<b_{\left(4,2\right)},\hbox{\ }
a_{\left(3^2\right)}<c_{\left(4,2\right)},\tag 4.1
$$
$$
b_{\left(3^2\right)}<a_{\left(4,2\right)},\hbox{\ }
b_{\left(3^2\right)}<b_{\left(4,2\right)},\hbox{\ }
b_{\left(3^2\right)}<c_{\left(4,2\right)},\tag 4.2
$$
$$
c_{\left(3^2\right)}<b_{\left(4,2\right)},\hbox{\ }
c_{\left(3^2\right)}<c_{\left(4,2\right)},\tag 4.3
$$
$$
d_{\left(3^2\right)}<c_{\left(4,2\right)},\tag 4.4
$$
and
$$
a_{\left(4,1^2\right)}<a_{\left(4,2\right)},\hbox{\ }
b_{\left(4,1^2\right)}<b_{\left(4,2\right)},\tag 4.5
$$
$$
c_{\left(4,1^2\right)}<b_{\left(4,2\right)},\hbox{\ }
d_{\left(4,1^2\right)}<c_{\left(4,2\right)},\hbox{\ }
e_{\left(4,1^2\right)}<c_{\left(4,2\right)}.\tag 4.6
$$

The inequalities (4.1) - (4.4) indicate the existence of the corresponding
(simple) substitution reactions among the $(4,2)$- and the $(3^2)$-derivatives,
and these substitution reactions can be used for complete identification of all
$(4,2)$-derivatives.  Indeed, two, three, and four $(3^2)$-products can be
synthesized from the $(4,2)$-derivatives which correspond to
$a_{\left(4,2\right)}$, $b_{\left(4,2\right)}$, and $c_{\left(4,2\right)}$,
respectively.

The following sets of structural formulae of $(3^2)$-derivatives can be
distinguished:
$$
\{a_{\left(3^2\right)}, b_{\left(3^2\right)}\},\hbox{\ }
\{c_{\left(3^2\right)}\},\hbox{\ }
\{d_{\left(3^2\right)}\}.
$$
Indeed, the products that correspond to the
elements of these sets can be synthesized from three, two, and one
$(4,2)$-derivatives, respectively.

The inequalities (4.5), (4.6) indicate the existence of the corresponding
(simple) substitution reactions among $(4,2)$- and $(4,1^2)$-derivatives, and
by means of these substitution reactions we can identify the following sets of
$(4,1^2)$-derivatives:
$$
\{a_{\left(4,1^2\right)}\},\hbox{\ }
\{b_{\left(4,1^2\right)}, c_{\left(4,1^2\right)}\},\hbox{\ }
\{d_{\left(4,1^2\right)}, e_{\left(4,1^2\right)}\}.
$$

Indeed, the product that corresponds to $a_{\left(4,1^2\right)}$ can be
synthesized only from the identifiable $a_{\left(4,2\right)}$, the products
that correspond to $b_{\left(4,2\right)}$ and $c_{\left(4,1^2\right)}$ can be
synthesized only from the identifiable $b_{\left(4,2\right)}$, and the products
that correspond to $d_{\left(4,2\right)}$ and $e_{\left(4,1^2\right)}$ can be
synthesized only from the identifiable $c_{\left(4,2\right)}$.

\heading
5. Genetic Relations: The Group $G$ has Order $6$ and is Dihedral

\endheading

In this section we describe the genetic relations of the molecule structure
under question when its Lunn-Senior's group $G$ of substitution isomerism
has order $6$, and is dihedral.  An instance is the molecule of cyclopropane
$C_3H_6$ (See [4]).  In accord with Theorem 2.1, we can suppose
$G=\langle (123)(456), (14)(26)(35)\rangle$.  Then
$T_{\left(4,2\right);G}=\{a_{\left(4,2\right)},b_{\left(4,2\right)},
c_{\left(4,2\right)},d_{\left(4,2\right)},\}$, where:

$a_{\left(4,2\right)}$ is the $G$-orbit
$$
\{(\{1,2,3,4\},\{5,6\}),
(\{1,2,3,5\},\{4,6\}),
(\{1,2,3,6\},\{4,5\}),
$$
$$
(\{2,4,5,6\},\{1,3\}),
(\{3,4,5,6\},\{1,2\}),
(\{1,4,5,6\},\{2,3\})\}
$$
of the tabloid
$A^{\left(4,2\right)}=
(\{1,2,3,4\},\{5,6\})$;

$b_{\left(4,2\right)}$ is the $G$-orbit
$$
\{(\{1,2,4,5\},\{3,6\}),
(\{2,3,5,6\},\{1,4\}),
(\{1,3,4,6\},\{2,5\})\}
$$
of the tabloid $B^{\left(4,2\right)}=(\{1,2,4,5\},\{3,6\})$;

$c_{\left(4,2\right)}$ is the $G$-orbit
$$
\{(\{1,2,4,6\},\{3,5\}),
(\{2,3,4,5\},\{1,6\}),
(\{1,3,5,6\},\{2,4\})\}
$$
of the tabloid
$C^{\left(4,2\right)}=(\{1,2,4,6\},\{3,5\})$;

$d_{\left(4,2\right)}$ is the $G$-orbit
$$
\{(\{1,2,5,6\},\{3,4\}),
(\{2,3,4,6\},\{1,5\}),
(\{1,3,4,5\},\{2,6\})\}
$$
of the tabloid
$D^{\left(4,2\right)}=(\{1,2,5,6\},\{3,4\})$.

We have

$T_{\left(3^2\right);G}=\{a_{\left(3^2\right)},b_{\left(3^2\right)},
c_{\left(3^2\right)},d_{\left(3^2\right)}\}$, where:

$a_{\left(3^2\right)}$ is the $G$-orbit
$$
\{(\{1,2,4\},\{3,5,6\}),
(\{2,3,5\},\{1,4,6\}),
(\{1,3,6\},\{2,4,5\}),
$$
$$
(\{2,4,5\},\{1,3,6\}),
(\{3,5,6\},\{1,2,4\}),
(\{1,4,6\},\{2,3,5\})\},
$$
of the tabloid
$A^{\left(3^2\right)}=
(\{1,2,4\},\{3,5,6\})$;

$b_{\left(3^2\right)}$ is the $G$-orbit
$$
\{(\{1,2,5\},\{3,4,6\}),
(\{2,3,6\},\{1,4,5\}),
(\{1,3,4\},\{2,5,6\}),
$$
$$
(\{1,4,5\},\{2,3,6\}),
(\{2,5,6\},\{1,3,4\}),
(\{3,4,6\},\{1,2,5\})\}
$$
of the tabloid
$B^{\left(3^2\right)}=
(\{1,2,5\},\{3,4,6\})$;

$c_{\left(3^2\right)}$ is the $G$-orbit
$$
\{(\{1,2,6\},\{3,4,5\}),
(\{2,3,4\},\{1,5,6\}),
(\{1,3,5\},\{2,4,6\}),
$$
$$
(\{3,4,5\},\{1,2,6\}),
(\{1,5,6\},\{2,3,4\}),
(\{2,4,6\},\{1,3,5\})\}
$$
of the tabloid
$C^{\left(3^2\right)}=
(\{1,2,6\},\{3,4,5\})$.

$d_{\left(3^2\right)}$ is the $G$-orbit
$$
\{(\{1,2,3\},\{4,5,6\}),
(\{4,5,6\},\{1,2,3\})\}
$$
of the tabloid
$D^{\left(3^2\right)}=
(\{1,2,3\},\{4,5,6\})$;

Moreover, we obtain

$T_{\left(4,1^2\right);G}=\{a_{\left(4,1^2\right)},b_{\left(4,1^2\right)},
c_{\left(4,1^2\right)},d_{\left(4,1^2\right)},e_{\left(4,1^2\right)}\}$, where:

$a_{\left(4,1^2\right)}$ is the $G$-orbit
$$
\{(\{1,2,3,4\},\{5\},\{6\}),
(\{1,2,3,5\},\{6\},\{4\}),
(\{1,2,3,6\},\{4\},\{5\}),
$$
$$
(\{2,4,5,6\},\{1\},\{3\}),
(\{3,4,5,6\},\{2\},\{1\}),
(\{1,4,5,6\},\{3\},\{2\})\}
$$
of the tabloid
$A^{\left(4,1^2\right)}=
(\{1,2,3,4\},\{5\},\{6\})$;

$b_{\left(4,1^2\right)}$ is the $G$-orbit
$$
\{(\{1,2,3,4\},\{6\},\{5\}),
(\{1,2,3,5\},\{4\},\{6\}),
(\{1,2,3,6\},\{5\},\{4\}),
$$
$$
(\{2,4,5,6\},\{3\},\{1\}),
(\{3,4,5,6\},\{1\},\{2\}),
(\{1,4,5,6\},\{2\},\{3\})\}
$$
of the tabloid
$B^{\left(4,1^2\right)}=
(\{1,2,3,4\},\{6\},\{5\})$;

$c_{\left(4,1^2\right)}$ is the $G$-orbit
$$
\{(\{1,2,4,5\},\{3\},\{6\}),
(\{2,3,5,6\},\{1\},\{4\}),
(\{1,3,4,6\},\{2\},\{5\}),
$$
$$
(\{1,2,4,5\},\{6\},\{3\}),
(\{2,3,5,6\},\{4\},\{1\}),
(\{1,3,4,6\},\{5\},\{2\})\}
$$
of the tabloid $C^{\left(4,1^2\right)}=(\{1,2,4,5\},\{3\},\{6\})$;

$d_{\left(4,1^2\right)}$ is the $G$-orbit
$$
\{(\{1,2,4,6\},\{3\},\{5\}),
(\{2,3,4,5\},\{1\},\{6\}),
(\{1,3,5,6\},\{2\},\{4\}),
$$
$$
(\{2,3,4,5\},\{6\},\{1\}),
(\{1,3,5,6\},\{4\},\{2\}),
(\{1,2,4,6\},\{5\},\{3\})\}
$$
of the tabloid
$D^{\left(4,1^2\right)}=
(\{1,2,4,6\},\{3\},\{5\})$;

$e_{\left(4,1^2\right)}$ is the $G$-orbit
$$
\{(\{1,2,5,6\},\{3\},\{4\}),
(\{2,3,4,6\},\{1\},\{5\}),
(\{1,3,4,5\},\{2\},\{6\}),
$$
$$
(\{1,3,4,5\},\{6\},\{2\}),
(\{1,2,5,6\},\{4\},\{3\}),
(\{2,3,4,6\},\{5\},\{1\})\}
$$
of the tabloid
$E^{\left(4,1^2\right)}=(\{1,2,5,6\},\{3\},\{4\})$.

This yields the inequalities
$$
A^{\left(3^2\right)}<A^{\left(4,2\right)}, \hbox{\ }
A^{\left(3^2\right)}<B^{\left(4,2\right)},\hbox{\ }
A^{\left(3^2\right)}<C^{\left(4,2\right)},\hbox{\ }
$$
$$
B^{\left(3^2\right)}<(123)(456)A^{\left(4,2\right)},\hbox{\ }
B^{\left(3^2\right)}<B^{\left(4,2\right)},\hbox{\ }
B^{\left(3^2\right)}<D^{\left(4,2\right)},\hbox{\ }
$$
$$
C^{\left(3^2\right)}<(132)(465)A^{\left(4,2\right)},\hbox{\ }
C^{\left(3^2\right)}<C^{\left(4,2\right)},\hbox{\ }
C^{\left(3^2\right)}<D^{\left(4,2\right)},
$$
$$
D^{\left(3^2\right)}<A^{\left(4,2\right)},
$$
and
$$
A^{\left(4,1^2\right)}<A^{\left(4,2\right)},\hbox{\ }
B^{\left(4,1^2\right)}<A^{\left(4,2\right)},\hbox{\ }
C^{\left(4,1^2\right)}<B^{\left(4,2\right)},
$$
$$
D^{\left(4,1^2\right)}<C^{\left(4,2\right)},\hbox{\ }
E^{\left(4,1^2\right)}<D^{\left(4,2\right)},
$$
so
$$
a_{\left(3^2\right)}<a_{\left(4,2\right)}, \hbox{\ }
a_{\left(3^2\right)}<b_{\left(4,2\right)},\hbox{\ }
a_{\left(3^2\right)}<c_{\left(4,2\right)},\tag 5.1
$$
$$
b_{\left(3^2\right)}<a_{\left(4,2\right)},\hbox{\ }
b_{\left(3^2\right)}<b_{\left(4,2\right)},\hbox{\ }
b_{\left(3^2\right)}<d_{\left(4,2\right)},\tag 5.2
$$
$$
c_{\left(3^2\right)}<a_{\left(4,2\right)},\hbox{\ }
c_{\left(3^2\right)}<c_{\left(4,2\right)},
c_{\left(3^2\right)}<d_{\left(4,2\right)},\tag 5.3
$$
$$
d_{\left(3^2\right)}<a_{\left(4,2\right)},\tag 5.4
$$
and
$$
a_{\left(4,1^2\right)}<a_{\left(4,2\right)},\hbox{\ }
b_{\left(4,1^2\right)}<a_{\left(4,2\right)},\tag 5.5
$$
$$
c_{\left(4,1^2\right)}<b_{\left(4,2\right)},\hbox{\ }
d_{\left(4,1^2\right)}<c_{\left(4,2\right)},\hbox{\ }
e_{\left(4,1^2\right)}<d_{\left(4,2\right)}.\tag 5.6
$$

The inequalities (5.1) - (5.4) indicate the existence of the corresponding
(simple) substitution reactions among the $(4,2)$- and the $(3^2)$-derivatives,
and the inequalities (5.5), (5.6) indicate the existence of the corresponding
(simple) substitution reactions among the $(4,2)$- and the
$(4,1^2)$-derivatives.

These substitution reactions can be used for distinguishing the products that
correspond to different sets from the following sets of structural formulae of
$(4,2)$-derivatives:  $$ \{a_{\left(4,2\right)}\},\hbox{\ }
\{b_{\left(4,2\right)},c_{\left(4,2\right)},d_{\left(4,2\right)}\}, $$ and from
the following sets of structural formulae of $(3^2)$-derivatives:  $$
\{a_{\left(3^2\right)}, b_{\left(3^2\right)},c_{\left(3^2\right)}\},\hbox{\ }
\{d_{\left(3^2\right)}\}.  $$ Indeed, it is enough to note that from the
product which corresponds to $a_{\left(4,2\right)}$ can be synthesized four
$(3^2)$-derivatives and from the products that correspond to the elements of
the set $\{b_{\left(4,2\right)}, c_{\left(4,2\right)},d_{\left(4,2\right)}\}$,
can be synthesized two $(3^2)$-derivatives.  The products that correspond to
the sets $\{a_{\left(3^2\right)}, b_{\left(3^2\right)},c_{\left(3^2\right)}\}$,
and $\{d_{\left(3^2\right)}\}$ can be synthesized from two and one
$(4,2)$-derivatives, respectively.

Using the above substitution reactions, we also can identify the products
corresponding to the following sets of structural
formulae of $(4,1^2)$-derivatives:
$$
\{a_{\left(4,1^2\right)},b_{\left(4,1^2\right)}\},\hbox{\ }
\{c_{\left(4,1^2\right)},d_{\left(4,1^2\right)},e_{\left(4,1^2\right)}\}.
$$
This is because both product that correspond to $a_{\left(4,1^2\right)}$ and
$b_{\left(4,1^2\right)}\}$ can be synthesized from the identifiable
$a_{\left(4,2\right)}$, and the products which correspond to
$c_{\left(4,1^2\right)}$, $d_{\left(4,1^2\right)}$, and
$e_{\left(4,1^2\right)}$ can be obtained from the products that correspond to
$b_{\left(4,2\right)}$, $c_{\left(4,2\right)}$, and $d_{\left(4,2\right)}$.

\heading {References}

\endheading

\noindent [1] W.  H\"asselbarth, The Inverse Problem of Isomer Enumeration, J.
Comput.  Chem.  8 (1987), 700 -- 717.

\noindent [2] V.  V.  Iliev, On Lunn-Senior's Mathematical Model of Isomerism
in Organic Chemistry.  Part I, MATCH - Commun.  Math.  Comput.  Chem.  40
(1999), 153 -- 186.

\noindent [3] V.  V.  Iliev, On the Inverse Problem of Isomer Enumeration:
Part I, Case of Ethane, MATCH - Commun.  Math.  Comput.  Chem.  43 (2001), 67
-- 77.

\noindent [4] V.  V.  Iliev, On the Inverse Problem of Isomer Enumeration:
Part II, Case of Cyclopropane, MATCH - Commun.  Math.  Comput.  Chem.  43
(2001), 79 -- 84.

\noindent [5] A. Kerber, Applied Finite Group Actions, Springer-Verlag, Berlin
1999.

\noindent [6] A.  C.  Lunn, J.  K.  Senior, Isomerism and Configuration, J.
Phys.  Chem.  33 (1929), 1027 - 1079.

\noindent [7] G.  P\'olya, Kombinatorische Anzahlbestimmungen f\"ur Gruppen,
Graphen und che\-mi\-sche Verbindungen, Acta Math.  68 (1937), 145 -- 254.
English translation:  G.  P\'olya and R.  C.  Read, Combinatorial Enumeration
of Groups, Graphs and Chemical Compounds, Springer-Verlag New York Inc., 1987.

\end